# Using Large Language Models as a Co-Author in Undergraduate Quantum Group Research

Jeffrey Kuan

Abstract

This article describes the use of Claude CLI and its Opus 4.6 model, as a tool for writing an entirely AI-generated mathematics research paper. The resulting paper is comparable in scope and quality to papers previously produced by advanced undergraduate students in eight-week summer REU programs advised by the author. The main result is a new explicit formula for a central element of $U_q(\mathfrak{so}_{12})$, which can be used for an interacting particle system with Markov duality. Using SageMath and a sparse PBW-basis pairing matrix that admits symbolic inversion, Claude reduced the central-element computation by several orders of magnitude: a calculation that took 60 hours in a 2023 Python implementation completed in under a minute on a laptop. The article reflects on the implications for undergraduate research mentorship: if generative AI can now produce research of REU caliber, advisors must select problems that better demonstrate the qualities valued by graduate admissions committees. Limitations including poor runtime estimates and literal handling of differing mathematical conventions are documented.

## Accessibility Statement

This PDF is WCAG2.1AA compliant and I tested it using VoiceOver. The PDF complies with Ohio Administrative Policy IT-09, Texas Administrative Code 206.70 (EFFECTIVE APRIL 18, 2020), Section 508 of the Rehabilitation Act and Title II of the Americans with Disabilities Act (effective April 24, 2027). For the technical details related to assistive technology usage, please visit my accessibility guidelines. Numbers in square brackets, such as [13], refer to entries in the bibliography at the end of the paper. For any issues, please send me an email.

## 1. Introduction

Over the last few years, many mathematicians have attempted to use generative AI for their research. The pace of progress has been striking: in the summer of 2025, AI models solved 5 of 6 International Mathematical Olympiad problems, and by February 2026, the "First Proof" challenge—in which 11 mathematicians posed 10 research–level lemmas spanning fields from algebraic combinatorics to representation theory—demonstrated AI solving at least 2 unambiguously and up to 8 partially [10]. Ernest Ryu used ChatGPT to prove a 42–year–old conjecture of Nesterov on optimization algorithms through an iterative three–day collaboration, in which approximately 80% of ChatGPT's suggestions were incorrect but several novel ideas led to the proof [11]. Separately, Ellenberg, Libedinsky, Plaza, Simental, and Williamson used Google DeepMind's AlphaEvolve to discover unexpected hypercube structures in Bruhat intervals of permutation groups—patterns that had gone unnoticed for 50 years [7]. As Daniel Litt remarked, "it's very likely that this technology is bigger than the computer" [18]. Terence Tao, who in October 2024 likened working with OpenAI's GPT-o1 to advising a "mediocre, but not completely incompetent, graduate student," declared in March 2026 that AI is now "ready for primetime" in mathematics because it "saves more time than it wastes" [26]. So far in the literature, generative AI has successfully been used to prove theorems, lemmas, and conjectures as part of longer research papers. Thus far, there has not been an entirely generative AI written research paper. This is in part because no reputable mathematics research journal nor arXiv accepts entirely generative AI research papers. In this article, I will describe my work in an entirely generative AI written research paper.

My primary motivation for doing this project comes from my undergraduate research advising. Over my career, I have advised over 20 students as part of seven summer REU programs. These programs have produced five research papers. As such, I have a responsibility to my students to assign problems that are more difficult than what generative AI can write. Therefore, I need to know what level of papers generative AI can produce.

The research paper was written with Claude CLI, after the release of Opus 4.6. The paper fits the level of what four advanced undergraduate students could produce in an 8 week summer program. The proof is computer assisted, using SAGE; there are a handful of journals which will accept computer assisted proofs, even though GenAI written papers are not yet accepted. Unlike some cutting–edge AI tools such as Google DeepMind's AlphaEvolve, which are only available to researchers with institutional partnerships [12], Claude CLI is publicly available, making this workflow accessible to any REU program.

The outline of this paper is as follows. Section 2 provides background on integrable probability and quantum groups, including the asymmetric simple exclusion process and the role of Casimir elements. Section 3 describes my undergraduate advising philosophy, the challenges posed by generative AI, and the impact of the undergraduate research on the wider community. Section 4 details the computation with Claude, including the use of SageMath and the troubleshooting process. Section 5 concludes with reflections on the future of AI–assisted mathematical research and its implications for undergraduate research mentorship.

## Acknowledgements

The author would like to gratefully thank The Ohio State University's Generative AI initiative, specifically OTDI (Office for Technology and Digital Innovation), ASC Tech (Arts and Sciences Technology); and The Ximera Foundation for providing a Claude Max subscription. The author acknowledges that some of this work was done on devices belonging to Texas A&M University.

## Conflict of Interest Disclosure

The namesake of the second author of the companion paper [16], "Tailor Swift Bot," is a company originally founded by the author while at Texas A&M University. Following the author's move to The Ohio State University, the company was restructured as a University

Technology Commercialization entity and renamed “Botting Accessibility With LaTeX” (B.A.W.L.) at the university’s request, in part due to not wanting to get sued by Taylor Swift. B.A.W.L. is a University Technology Commercialization Company (UTCC) as defined in Chapter 3335-13 of The Ohio State University’s bylaws and rules. The views and products of B.A.W.L. are those of the company and do not necessarily reflect the views of, nor are they endorsed or approved by, The Ohio State University.

## 2. Research Area: Integrable Probability and Quantum Groups

### 2.1 Integrable Probability Background

Integrable probability is a branch of probability theory that studies models which have exact solutions. This is in contrast to more general probability theory, where most models usually require estimates. The term “integrable” has its origins in classical Hamiltonian mechanics, where an integrable system is one possessing enough conserved quantities to solve the equations of motion explicitly. A (somewhat) simple example is the harmonic oscillator: an object on a frictionless spring, whose total energy is conserved and whose position can be written in closed form. In contrast, the three-body problem is not integrable, and its solutions are notoriously chaotic (it is well known enough to have a novel and a TV series named after it).

During the 1980s, the “Soviet school” of mathematical physics introduced quantum mechanics into integrable systems. In this context, the notion of “conserved quantities” is replaced by “commuting operators.” These operators no longer commute, a concept most commonly understood by the idea that you cannot simultaneously measure both position and momentum. The most famous mathematician from this school is Vladimir Drinfel’d, who won a Fields Medal in part for this work. Since then, the “Drinfel’d–Jimbo” quantum groups have been generalized in the algebraic framework of Hopf algebras.

In the context of probability, one of the canonical models is the asymmetric simple exclusion process (ASEP), introduced by Spitzer [23] in 1970, and earlier in biology in 1968 [20]. The word “exclusion” means that at most one particle can occupy a site. The word “simple” means that particles jump at most distance one at a time. “Asymmetric” means that the jump rates are asymmetric in spatial directions. The nontrivial interaction occurs when a particle attempts to jump to an occupied site; in such cases the jumps are blocked.

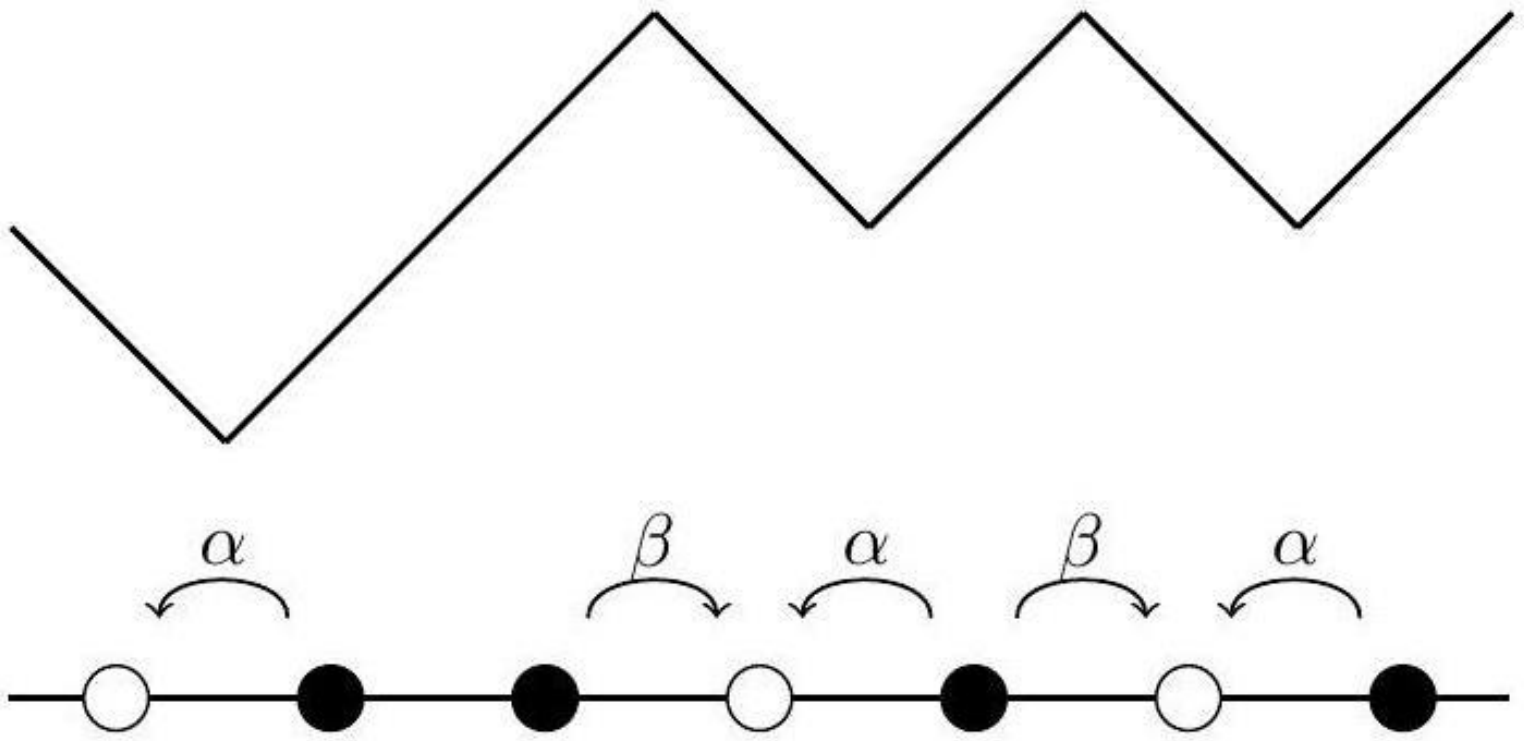

*Caption: A configuration of the asymmetric simple exclusion process (ASEP) on a one-dimensional lattice, with the corresponding height function.*

The jump rates can be collected in a square matrix called the *generator*. The dimension of the generator is the number of possible states, which for ASEP on $L$ lattice sites is $2^L$. The generator $A$ is a matrix whose rows sum to 0 and has non–negative off-diagonal entries. The reason the rows sum to 0 is that the transition matrix at time $t$ is defined by $e^{tA} = 1 + tA + t^2A^2/2! + \cdots$, which has rows which sum to 1. This corresponds to probabilities summing to 1. The off-diagonal entries need to be non–negative because probabilities are always non–negative.

After Drinfeld's work, it was discovered in the early 90's that ASEP is related to the Heisenberg XXZ model (which models magnets) and the quantum group $U_q(\mathfrak{sl}_2)$. Starting in the late 90's and early 2000's, it was discovered that integrable models display a phenomenon called **universality**. Even in popular culture, the concept of universality is well–known enough that there is a meme based on it in the case of the Bell curve and the normal distribution.

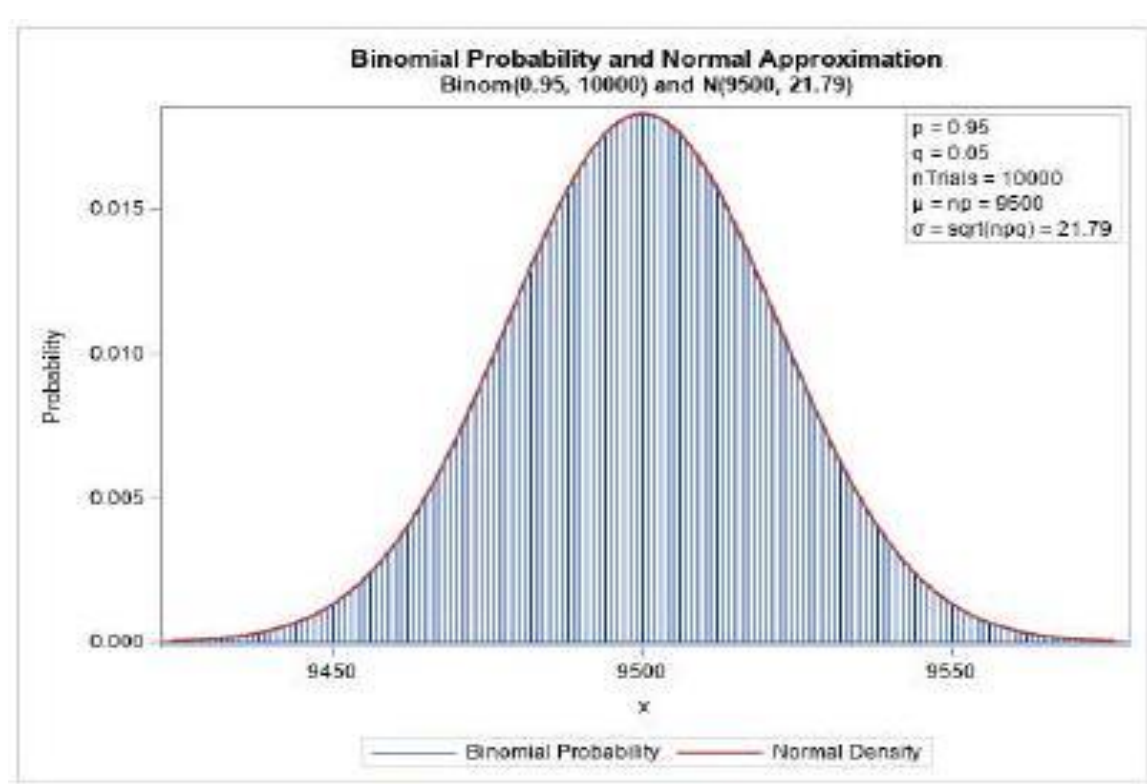

*Caption: An overlay of the bell curve (the standard Normal distribution) and the binomial distribution, illustrating universality.*

In the setting of ASEP, the universal distribution is the Tracy–Widom distribution [24, 25], which has 1/3 scaling exponent instead of 1/2. The two images illustrate the two universal distributions.

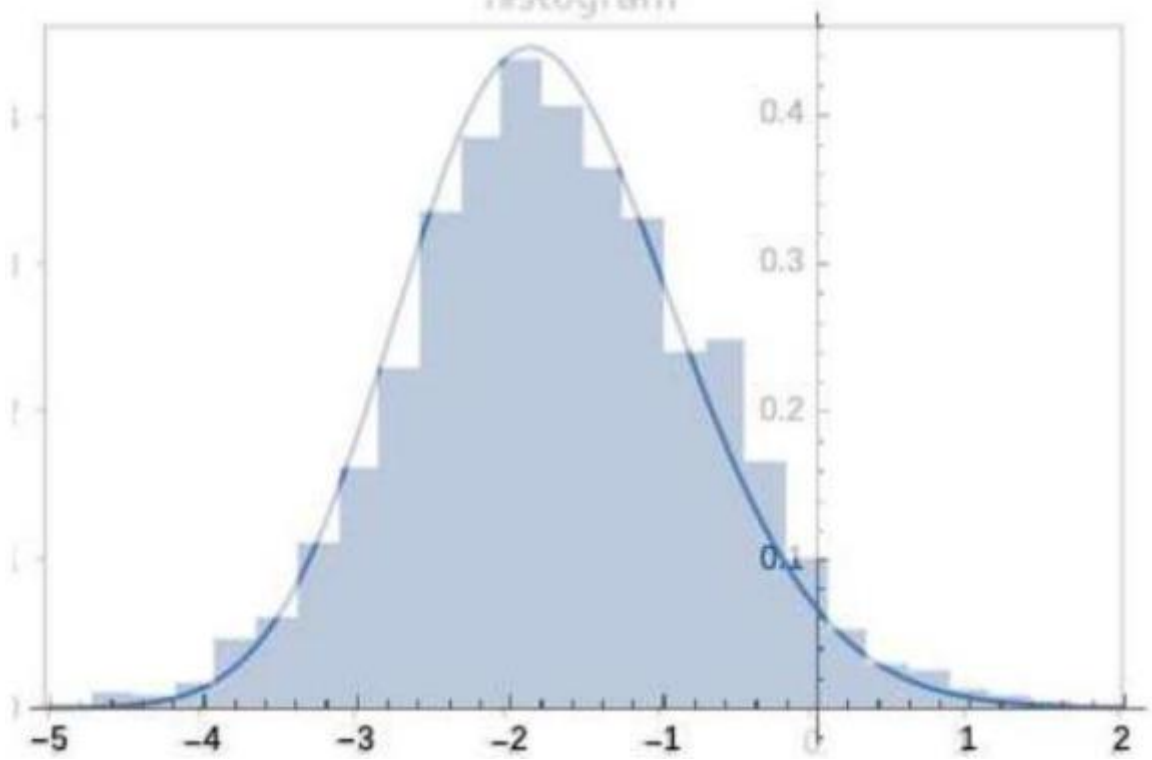

*Caption: The Tracy–Widom distribution. It resembles a skewed bell curve, with a longer tail on the left side and a sharper dropoff on the right.*

The proof of this result uses the exact solutions of ASEP, which can be written using contour integrals. A later proof [3] was found using a probabilistic method called Markov duality. Despite its ubiquity, Markov duality has been called a “black art.” A large part of my research has been turning this art into a science, by leveraging the underlying algebraic machinery to prove Markov duality.

## 2.2 Algebraic Background

The quantum group $U_q(\mathfrak{sl}_2)$ can be generalized to more complex algebraic objects. The notation $\mathfrak{sl}_2$ refers to the Lie algebra of $2 \times 2$ matrices with trace zero. It is an example of a finite–dimensional simple Lie algebra. The set of all finite–dimensional simple Lie algebras can be categorized through something called Dynkin diagrams.

$A_n$ $E_6$

$B_n$

$E_7$

$C_n$

$E_8$

$D_n$

$F_4$

$G_2$

*Alt Text: Dynkin diagrams for the four infinite families of finite–dimensional simple Lie algebras: $A_n$, $B_n$, $C_n$, and $D_n$, with ellipses indicating the variable number of vertices. The $A_n$ graph is a straight line, while the $D_n$ graph branches at the end. The exceptional Lie algebras $E_6$, $E_7$, $E_8$ are similar to $D_n$ in that they branch, but the branch is in the middle rather than at the end. The remaining exceptional Lie algebras $F_4$ and $G_2$ are straight lines with 4 and 2 vertices, respectively, with a double or triple edge.*

There are four infinite series as well as five “exceptional” Lie algebras. For the case of type $A_n$ Lie algebras, Belitsky–Schütz [1] and I [13] proved a Markov duality for the so–called “multi–species ASEP” introduced by Liggett in 1976 [17].

The calculation involves an “explicit” expression for a certain algebraic element called the *Casimir element*. Roughly speaking, an $m \times m$ matrix algebra acts on an $m$–dimensional vector space by usual multiplication, which can be generalized to representations. In the case of $\mathfrak{sl}_2$ and ASEP, the representation of any element of $\mathfrak{sl}_2$ on the $2^L$–dimensional vector space $V \otimes V \otimes \cdots \otimes V$ ($L$ tensor products) is a $2^L \times 2^L$ matrix. The matrix size corresponds to the $2^L$ possible states of ASEP on $L$ lattice sites. This matrix has a particular name, called the *Hamiltonian*. If the element is the Casimir element, then a change of basis results in the ASEP generator.

What constitutes an “explicit” expression will depend on one’s background and perspective, and it is actually due to these differing perspectives that I am able to assign open problems to students. Algebraicists prefer proving a statement that holds for all finite simple Lie algebras (i.e. all Lie algebras expressed in Figure 1), and such formulas for the Casimir do exist (such as in Drinfeld’s R-matrix construction [6]). However, from a probabilistic perspective, one needs to expand the Casimir over monomials, so one can calculate the matrix entries and check that the probabilities are nonnegative. Because of these differing perspectives, explicit (to a probabilist) expressions for the Casimir of other simple Lie algebras are not in the literature.

In fact, there is no reason to expect the off-diagonal entries to be nonnegative. This is not an insurmountable issue. If a row has negative off-diagonal entries, we view this row to be “problematic”. We then adjust the change of basis matrix so that the probability of entering a problematic state is 0. In this sense, problematic states are discarded.

## 3. Undergraduate Advising

### 3.1 Undergraduate Advising Philosophy

Throughout my career, I have advised seven summer REU projects: 2016 and 2017 (while a NSF postdoc at Columbia University) and 2020–2024 (while tenure track at Texas A&M University). Each summer project was with a group of three or four students; collectively, they have produced five papers, which have all been accepted to research journals. Through experience, I’ve determined the projects work best when some students are interested in algebra and some students are interested in probability, as it fosters productivity in a group setting. Because collaboration and collegiality is essential to the success of the program, I heavily value the students’ personal statements when evaluating applications.

Arguably the most important part of advising undergraduates is to help them develop a strong graduate school application. In the past, REU’s were meant more for research *experience*, although nowadays there is more of an *expectation* that students have research publications. Therefore, the “output” of the REU project should demonstrate qualities which are valued by graduate school application committees. Some obvious qualities are creativity and having a strong mathematical background. A slightly less obvious quality is perseverance, since maintaining interest throughout an 8 week program is much different than working on weekly homework assignments.

Historically, a large proportion of my students have been successful in their graduate school applications. Some students chose not to apply to graduate school in math, but every student who applied received at least one acceptance. While each summer project focuses on different problems, depending on the interests of the students, the presence of a research paper as output has resulted in graduate school acceptances. In 2020 and 2023, the students used a Python assisted proof to calculate the Casimir element. At the time, generative AI had not yet developed enough that the use of a Python script significantly impacted their success in graduate school applications.

The students already knew Python but did not know SageMath. For the smaller calculations in the 2020 and 2023 projects, SageMath would have sped up the computation but was strictly speaking not necessary. Additionally, I do not know SageMath, so I could not help the students troubleshoot had they chosen to use it.

With recent advances in Generative AI, especially with "vibe coding", one can question whether computer assisted proofs truly demonstrate the qualities desired by a graduate school admissions committee. As Terence Tao has observed, current AI is "very good at scouring big lists of problems for low–hanging fruit" [12]. AI can now instantly solve many homework–level problems, undermining the traditional role of problem sets in developing mathematical maturity. As Aravind Asok of USC has noted, he stopped assigning homework because graduate students use AI to generate answers [21]. If this is already happening at the graduate level, it raises questions for undergraduate research advisors. In my role as a research advisor, I have a responsibility to position my students to be successful in their careers, which means I need to assign problems that can illuminate their talents and abilities, by assigning problems that are not low–hanging fruit. As Tao has argued, "the scarce skill becomes choosing the right problem, designing the workflow, and checking the result" [26]. Choosing the right problem is a difficult task for undergraduate research advisors. My track record of five published papers and a 100% graduate school acceptance rate among applicants suggests that problem selection has been a strength of my advising, and I intend to preserve that strength by adapting to the new landscape. Empirical evidence supports the importance of mentorship in this context: a 2026 meta–analysis of generative AI in mathematics education found that students with teacher support showed 18 times larger learning gains than those working with AI alone [19].

## 3.2 Impact on the Wider Research Community

A "secondary" value of REUs is the actual paper itself. In this case, the papers were published in research journals that were not undergraduate research journals. The type D ASEP consists of particles of two types, similar to how the type $A_n$ ASEP is the multispecies ASEP introduced by Liggett in 1976 [17]. However, in type D ASEP the parameter $n$ actually affects the speed and not the state space. This is due to discarding “problematic” particle configurations. The duality function and reversible measures for the type D ASEP are products of two copies of the duality function and reversible measures for usual ASEP. This suggests a convergence to a product of two independent copies of Tracy–Widom, which would be mildly interesting to the integrable probability community at large. A natural question that has come up about my research is the strength of the interaction between the different species of particles. In the trivial case of independence (no interactions), it is straightforward that the asymptotics are two independent Tracy–Widom distributions. This model gives a qualitative description of the minimal strength of interaction needed for nontrivial asymptotics.

From an algebraic perspective, it is also worth asking how $n$ affects the interacting particle system for $n = 6,7,8$ in the exceptional Lie algebras $E_6$, $E_7$, $E_8$. However, the dimensions of these algebras are much larger, and much less intuitive as matrix algebras. The Lie algebra $D_6$ is a natural “bridge” between type $D$ and type $E$.

The next section will describe new mathematical research, written using Claude CLI after the release of Opus 4.6. The results and proof are comparable to previously published REU papers [15, 14, 2, 22, 4], of which [14, 22] use computer assistance. Based on this experience, I conclude that my previously assigned problems are no longer appropriate in the context of summer REUs.

# 4. The Computation with Claude

## 4.1 Comparison with Previous Methods

The 2023 REU paper has Python code available on GitHub, which completed in about 45 minutes on a Mac Mini for $U_q(\mathfrak{so}_{10})$. The first step was to test the code for $U_q(\mathfrak{so}_{12})$. As a test case, the code first attempted to find a basis of a weight space of $U_q(\mathfrak{so}_{12})$. The script took 60 hours to complete on an Amazon Large EC2 instance, and failed to match the correct dimension. This was due to floating point errors; the script calculates the inverse of large matrices at various values of $q$, then uses interpolation to determine the precise formulas. The script only found 155 basis elements in a 200–dimensional weight space. Clearly, a new method was needed.

## 4.2 Claude and SageMath

One method, suggested to me by Scott Blair during the 2025 Ross Program, is to use SageMath, and specifically QuaGroup, for the higher–dimensional Lie algebras. Historically, this would’ve required me to learn SageMath; however, for the purposes of this project, I intentionally chose not to learn anything about SageMath, to test out if Claude CLI could generate a paper with minimal knowledge.

The result produced by Claude was significantly more efficient. For example, the aforementioned 60 hour computation took less than one minute on my MacBook Air. The bottleneck in the previous calculation was inverting a large matrix. In the new calculation, the pairing matrix $M$ is diagonal. The remaining "change of basis" matrix $B$ is very sparse [16]. For the largest weight space (dimension 275), the $275 \times 275$ matrix that we invert has only 2267 nonzero entries out of 75,625 total, about 3% nonzero. This made the computation much faster and did not require interpolation.

The key insight is that the pairing matrix, calculated on the Sage PBW basis elements, yields a diagonal matrix. I intentionally chose not to read the Sage or QuaGroup documentation to understand why. Every computed case resulted in a diagonal matrix and I assume there is some deeper mathematical reason.

With this insight, one can invert a sparser matrix. We expand monomials (of algebra generators) over the PBW basis. This "change of basis" matrix is sparse and we invert this matrix instead. This calculation can be done symbolically and does not require interpolation, as in the previous script from 2023.

One way to verify that an algebraic element is central is by checking that it commutes with all algebra elements. This can be done by hand (in theory). In this case, I asked Claude to write a Sage script to verify that the commutation holds. There is also a more "human" intuition for the correctness of the result. Namely, the coefficients were all relatively simple functions of $q$, with numerical coefficients being somewhat small integers. An erroneous script would more likely result in very "complicated" formulas.

There is, of course, a (probably true) stereotype that nobody actually reads papers, including the referees. Arguably, having Claude verify a computation is more trustworthy than having a human verify a calculation. This parallels a broader trend: as Tao notes in the Quanta Magazine article [12], "AI without validation is too unreliable," and human verification of proofs remains essential. As Emily Riehl of Johns Hopkins has observed, "there's nothing in the way that large language models are designed" to ensure mathematical correctness [21]. Tao and collaborators found a similar lesson when using AlphaEvolve on 67 mathematical problems: without carefully designed, exploit–proof verification code, the AI would find ways to game the scoring function rather than produce genuine solutions [8]. Another possible source of verification is Lean, a proof assistant and functional programming language. As an example of its usage, it is currently being used by Kevin Buzzard and over 60 contributors to formalize Fermat's Last Theorem [21]. As Buzzard has noted, "solving hallucinations is not something that's going to happen anytime soon" [5]. Indeed, recent theoretical work has proven that hallucinations are mathematically inevitable in large language models [27]. Lean may be a topic pursued in a future paper.

### 4.3 Troubleshooting

As a human analogy, an advisor should decide on how often to meet with a student (or postdoc) to check in on their current progress. A more independent researcher may need fewer meetings and can reach out as needed, while a more novice researcher would need more frequent meetings to ensure they are on the right path.

#### *4.3.1 Poor Time Estimates*

A recurring source of frustration was Claude's inability to accurately estimate computation time. On one notable occasion, Claude predicted that a computation would take "microseconds per entry" or complete "instantly," when the actual runtime was minutes or hours. In one session, after several such predictions proved wildly optimistic, I noted that the estimate of "microseconds" was completely wrong.

This is a limitation for an AI assistant. When planning a workflow, accurate time estimates are essential for deciding whether to run a computation interactively or as a batch job, and whether intermediate results should be saved to disk. Claude's estimates were based on poor guesses and were presented with unwarranted confidence. After this experience, I asked Claude to run a small unit test and extrapolate from that data, or explicitly say "I don't know how long this will take."

#### *4.3.2 Same versus Different*

Often in colloquial mathematics conversations, two technically different objects are understood to be "morally" the same. For example, in the asymmetric simple exclusion process, some authors take the drift to the left and others take drift to the right. There are also analogously "different" conventions for the coproduct of the quantum groups. Human beings understand these to be conventions and the same in essence. However, Claude CLI took definitions and conventions more literally than most humans would. This resulted in several problems that needed troubleshooting. As Timothy Gowers has noted, "humans seem to be able to find proofs with a remarkably small amount of backtracking," and he has expressed dissatisfaction with programs that find proofs "after carrying out large searches, even if those searches are small enough to be feasible" [9]. In other words, humans are efficient at diagnosing *why* an approach fails and adjusting accordingly, rather than starting over from scratch. The root vector convention bug described below is a concrete example: Claude could not diagnose why the factor of $-1/q$ appeared, but a human could trace it to differing braid group conventions.

For example, our root vector expansions used $T_1(F_2) = F_1F_2 - q^{-1}F_2F_1$, while Sage uses $T_1(F_2) = F_2F_1 - qF_1F_2$. These differ by a factor of $-1/q$, although a human would recognize those as essentially the same. The fix was to build root vector matrices directly from Sage's PBW commutation relations, although Claude had a difficult time identifying the problem.

## 5. Conclusion

This paper described the use of Claude CLI, using Opus 4.6, as a digital tool in writing an undergraduate research-level mathematics paper. The resulting paper [16] is comparable to papers previously produced by advanced undergraduate students in 8-week summer REU programs advised by the author. Claude streamlined the computer calculations for $U_q(\mathfrak{so}_{12})$ central elements by several orders of magnitude.

The implications for undergraduate research advising are significant. If a generative AI tool can produce research comparable to what REU students have published, then the problems I have previously assigned no longer adequately demonstrate the qualities that graduate school admissions committees seek—creativity, perseverance, and deep mathematical understanding. This does not mean that undergraduate research is obsolete; rather, it means that advisors must be more deliberate in selecting problems that require the personal qualities valued by the mathematical community.

The troubleshooting experience also reveals important limitations of AI as a mathematical tool. Claude's poor time estimates, its literal interpretation of differing conventions, and its inability to distinguish "morally equivalent" mathematical objects all required human intervention to resolve. These are not bugs that will be fixed with larger models; they reflect fundamental differences between how humans and language models process mathematical knowledge. A human mathematician carries implicit understanding of conventions, approximations, and the relative importance of different details. Claude does not.

Looking forward, several directions present themselves. Formal verification using Lean could provide a more rigorous check on AI–generated proofs than the computational verification used here. The exceptional Lie algebras $E_6, E_7, E_8$ offer natural extensions of the type D computation, though their larger dimensions present substantial challenges even for AI–assisted methods. Most importantly, the rapidly evolving capabilities of generative AI mean that the boundary between "problems AI can solve" and "problems that require human insight" will continue to shift, requiring ongoing vigilance from research mentors.

The central question for mathematics education in this era is how to ensure that students develop genuine mathematical understanding alongside their facility with generative AI tools.